\documentclass[12pt,twoside]{article} 
\usepackage[english]{babel}
\usepackage[latin1]{inputenc} 
\usepackage{amsmath}
\usepackage{amssymb,amsfonts}
\usepackage{graphicx}                   
\usepackage{times,amssymb,amscd} 

\newcommand{\bA}{\mathbf{A}}

\newcommand{\bE}{\mathbf{E}}

\newcommand{\bH}{\mathbf{H}}

\newcommand{\bR}{\mathbf{R}}
\newcommand{\bS}{\mathbf{S}}
\newcommand{\bV}{\mathbf{V}}

\newcommand{\be}{\mathbf{e}}

\newcommand{\br}{\mathbf{r}}
\newcommand{\bx}{\mathbf{x}}
\newcommand{\by}{\mathbf{y}}

\newcommand{\BV}{\boldsymbol{V}}

\newcommand{\Be}{\boldsymbol{e}}
\newcommand{\Bu}{\boldsymbol{u}}
\newcommand{\Bv}{\boldsymbol{v}}

\newcommand{\cP}{\mathcal{P}}
\newcommand{\cS}{\mathcal{S}}

\newcommand{\EUC}{\bE^3}

\newcommand{\SXR}{\bS^2\!\times\!\bR}
\newcommand{\HXR}{\bH^2\!\times\!\bR}

\newcommand{\NIL}{\mathbf{Nil}}
\newcommand{\SOL}{\mathbf{Sol}}

\newtheorem{Definition}{Definition}[section]
\newtheorem{Remark}{Remark}[section]

\begin{document}
\pagestyle{myheadings}
\markboth{\centerline{Jen\H o Szirmai}}
{On lattice coverings of $\NIL$ space by congruent geodesic balls}
\title
{On lattice coverings of $\NIL$ space by congruent geodesic balls }

\author{Jen\H o Szirmai \\
\normalsize Budapest University of Technology and \\ 
\normalsize Economics Institute of Mathematics, \\
\normalsize Department of Geometry \\
\normalsize Budapest, P. O. Box: 91, H-1521 \\
\normalsize szirmai@math.bme.hu
\date{\normalsize{\today}}}

%%%%%%%%%%%%%%%%%%%%%%%%%%%%%%%%%%%%%%%%%%%%
%%AMS Classification 2000
%% The 1-st classification is obligatory, the 2-nd classification is optional
%% \subjclass{primary}{secondary}      f.e. \subjclass{35R35, 49N50}{}
%%{52C17, 52C22}{}
%%%%%%%%%%%%%%%%%%%%%%%%%%%%%%%%%%%%%%%%%%%%

\maketitle
\begin{abstract}

The $\NIL$ geometry, which is one of the eight 3-dimensional  
Thurston geometries, can be derived from {W. Heisenberg}'s famous real matrix group.

The aim of this paper to study {\it lattice coverings} in $\NIL$ space. We introduce the notion of the density of considered coverings and
give upper and lower estimations to it, moreover we formulate a conjecture for the ball arrangement of the least dense lattice-like geodesic
ball covering and give its covering density $\Delta\approx 1.42900615$. 

The homogeneous 3-spaces have a unified interpretation in the projective 3-sphere and 
in our work we will use this projective model of the $\NIL$ geometry. 

\end{abstract}
%%%%%%%%%%%%%%%%%%%%%%%%%%%%%%%%%%%%%%%%%%%%

%%%%%%%%%%%%%%%%%%%%%%%%%%%%%%%%%%%%%%%%%%% 
\newtheorem{theorem}{Theorem}[section]
\newtheorem{corollary}[theorem]{Corollary}
\newtheorem{conjecture}[theorem]{Conjecture}
\newtheorem{lemma}[theorem]{Lemma}
\newtheorem{exmple}[theorem]{Example}
\newtheorem{defn}[theorem]{Definition}
\newtheorem{rmrk}[theorem]{Remark}
\newenvironment{definition}{\begin{defn}\normalfont}{\end{defn}}
\newenvironment{remark}{\begin{rmrk}\normalfont}{\end{rmrk}}
\newenvironment{example}{\begin{exmple}\normalfont}{\end{exmple}}
 \newenvironment{acknowledgement}{Acknowledgement}

%%%%%%%%%%%%%%%%%%%%%%%%%%%%%%%%%%%%%%%%%%%%%%%%%%%%%%%%%%%%%%%%%%%% 

%============================================================================%
%                             the main article                               %
%============================================================================%

%%%%%%%%%%%%%%%%%%%%%%%%%%%%%%%%%%%%%%%%%%%%%%%%%%%%%%%%%%%%%%%%%%%%%%%%%%%%%%
\section{Notions of the Nil geometry}

In this Section we summarize the significant notions and denotations of the $\NIL$ geometry (see \cite{M97}, \cite{Sz07-2}).

The $\NIL$ geometry is a homogeneous 3-space derived from the famous real matrix group $\mathbf{L(\mathbb{R})}$ discovered by {Werner Heisenberg}. 
The Lie theory with the methode of the projective geometry makes possible to investigate and to describe this topic.

The left (row-column) multiplication of Heisenberg matrices
     \begin{equation}
     \begin{gathered}
     \begin{pmatrix}
         1&x&z \\
         0&1&y \\
         0&0&1 \\
       \end{pmatrix}
       \begin{pmatrix}
         1&a&c \\
         0&1&b \\
         0&0&1 \\
       \end{pmatrix}
       =\begin{pmatrix}
         1&a+x&c+xb+z \\
         0&1&b+y \\
         0&0&1 \\
       \end{pmatrix}
      \end{gathered} \tag{1.1}
     \end{equation}
defines "translations" $\mathbf{L}({\mathbb{R}})= \{(x,y,z): x,~y,~z\in \mathbb{R} \}$ 
on the points of the space $\NIL= \{(a,b,c):a,~b,~c \in \mathbb{R}\}$. 
These translations are not commutative in general. The matrices $\mathbf{K}(z) \vartriangleleft \mathbf{L}$ of the form
     \begin{equation}
     \begin{gathered}
       \mathbf{K}(z) \ni
       \begin{pmatrix}
         1&0&z \\
         0&1&0 \\
         0&0&1 \\
       \end{pmatrix}
       \mapsto (0,0,z)  
      \end{gathered}\tag{1.2}
     \end{equation} 
constitute the one parametric centre, i.e. each of its elements commutes with all elements of $\mathbf{L}$. 
The elements of $\mathbf{K}$ are called {\it fibre translations}. $\NIL$ geometry of the Heisenberg group can be projectively 
(affinely) interpreted by the "right translations" 
on points as the matrix formula 
     \begin{equation}
     \begin{gathered}
       (1;a,b,c) \to (1;a,b,c)
       \begin{pmatrix}
         1&x&y&z \\
         0&1&0&0 \\
         0&0&1&x \\
         0&0&0&1 \\
       \end{pmatrix}
       =(1;x+a,y+b,z+bx+c) 
      \end{gathered} \tag{1.3}
     \end{equation} 
shows (see (1.1)). Here we consider $\mathbf{L}$ as projective collineation group with right actions in homogeneous coordinates.
We will use the Cartesian homogeneous coordinate simplex $E_0(\be_0)$$E_1^{\infty}(\be_1)$$E_2^{\infty}(\be_2)$
$E_3^{\infty}(\be_3), \ (\{\be_i\}\subset \bV^4$ \ with the unit point $E(\be = \be_0 + \be_1 + \be_2 + \be_3 ))$ 
which is distinguished by an origin $E_0$ and by the ideal points of coordinate axes, respectively. 
Moreover, $\by=c\bx$ with $0<c\in \mathbb{R}$ (or $c\in\mathbb{R}\setminus\{0\})$
defines a point $(\bx)=(\by)$ of the projective 3-sphere $\cP \cS^3$ (or that of the projective space $\cP^3$ where opposite rays
$(\bx)$ and $(-\bx)$ are identified). 
The dual system $\{(\Be^i)\}, \ (\{\be^i\}\subset \BV_4)$ describes the simplex planes, especially the plane at infinity 
$(\Be^0)=E_1^{\infty}E_2^{\infty}E_3^{\infty}$, and generally, $\Bv=\Bu\frac{1}{c}$ defines a plane $(\Bu)=(\Bv)$ of $\cP \cS^3$
(or that of $\cP^3$). Thus $0=\bx\Bu=\by\Bv$ defines the incidence of point $(\bx)=(\by)$ and plane
$(\Bu)=(\Bv)$, as $(\bx) \text{I} (\Bu)$ also denotes it. Thus {$\NIL$} can be visualized in the affine 3-space $\bA^3$
(so in $\bE^3$) as well.

The translation group $\mathbf{L}$ defined by formula (1.3) can be extended to a larger group $\mathbf{G}$ of collineations,
preserving the fibering, that will be equivalent to the (orientation preserving) isometry group of $\NIL$. 

In \cite{M06} E.~Moln\'ar has shown that 
a rotation trough angle $\omega$
about the $z$-axis at the origin, as isometry of $\NIL$, keeping invariant the Riemann
metric everywhere, will be a quadratic mapping in $x,y$ to $z$-image $\overline{z}$ as follows:
     \begin{equation}
     \begin{gathered}
       \mathcal{R}=\br(O,\omega):(1;x,y,z) \to (1;\overline{x},\overline{y},\overline{z}); \\ 
       \overline{x}=x\cos{\omega}-y\sin{\omega}, \ \ \overline{y}=x\sin{\omega}+y\cos{\omega}, \\
       \overline{z}=z-\frac{1}{2}xy+\frac{1}{4}(x^2-y^2)\sin{2\omega}+\frac{1}{2}xy\cos{2\omega}.
      \end{gathered} \tag{1.4}
     \end{equation}
This rotation formula $\mathcal{R}$, however, is conjugate by the quadratic mapping $\mathcal{M}$ to the linear rotation $\Omega$ in (1.5) as follows
     \begin{equation}
     \begin{gathered}
       \mathcal{M}: \ \ (1;x,y,z) \stackrel{\mathcal{M}}{\longrightarrow} (1; x',y',z')=(1;x,y,z-\frac{1}{2}xy) \ \ \text{to} \\
       \Omega: \ \ (1;x',y',z') \stackrel{\Omega}{\longrightarrow} (1;x",y",z")=(1;x',y',z')
       \begin{pmatrix}
         1&0&0&0 \\
         0&\cos{\omega}&\sin{\omega}&0 \\
         0&-\sin{\omega}&\cos{\omega}&0 \\
         0&0&0&1 \\
       \end{pmatrix}, \\
       \text{with} \ \ \mathcal{M}^{-1}: (1;x",y",z") \stackrel{\mathcal{M}^{-1}}{\longrightarrow}  (1; \overline{x}, \overline{y},\overline{z})=(1; x",y",z"+\frac{1}{2}x"y").
      \end{gathered} \tag{1.5}
     \end{equation}
This quadratic conjugacy modifies the $\NIL$ translations in (1.1), as well. Now a translation with $(X,Y,Z)$ in (1.3) instead of $(x,y,z)$ will be changed 
by the above conjugacy to the translation 
     \begin{equation}
     \begin{gathered}
        (1;x,y,z) \longrightarrow (1; \overline{x}, \overline{y},\overline{z})=(1; x,y,z)
        \begin{pmatrix}
         1&X&Y&Z-\frac{1}{2}XY \\
         0&1&0&-\frac{1}{2}Y \\
         0&0&1&\frac{1}{2}Y \\
         0&0&0&1 \\
       \end{pmatrix}, \\
             \end{gathered} \tag{1.6}
     \end{equation}
that is again an affine collineation.
We shall use the following important classification theorem.
\begin{theorem}[E.~Moln\'ar \cite{M06}]

{\bf (1)} Any group of $\NIL$ isometries, containing a 3-dimensional translation lattice, 
is conjugate by the quadratic mapping in (1.5) to an affine group of the affine (or Euclidean) space $\bA^3=\EUC$ 
whose projection onto the (x,y) plane is an isometry group of $\bE^2$. Such an affine group preserves a plane
$\to$ point polarity of signature $(0,0,\pm 0,+)$.
\newline
{\bf (2)} Of course, the involutive line reflection about the $y$ axis
     \begin{equation}
     \begin{gathered}
       (1;x,y,z) \to (1;-x,y,-z), 
      \end{gathered} \notag
     \end{equation}
preserving the Riemann metric, and its conjugates by the above isometries in {(\bf 1)} (those of the identity component)
are also {$\NIL$}-isometries. Orientation reversing {$\NIL$}-isometry does not exist.
\end{theorem}

The geodesic curves of the $\NIL$ geometry are generally defined as having locally minimal arc length between their any two (near enough) points. 
The equation systems of the parametrized geodesic curves $g(x(t),y(t),z(t))$  in our model can be determined by the 
general theory of Riemann geometry:
We can assume, that the starting point of a geodesic curve is the origin because we can transform a curve into an 
arbitrary starting point by translation (1.1); 
\begin{equation}
\begin{gathered}
        x(0)=y(0)=z(0)=0; \ \ \dot{x}(0)=c \cos{\alpha}, \ \dot{y}(0)=c \sin{\alpha}, \\ \dot{z}(0)=w; \ - \pi \leq \alpha \leq \pi. \notag
\end{gathered}
\end{equation}
The arc length parameter $s$ is introduced by
\begin{equation}
 s=\sqrt{c^2+w^2} \cdot t, \ \text{where} \ w=\sin{\theta}, \ c=\cos{\theta}, \ -\frac{\pi}{2}\le \theta \le \frac{\pi}{2}, \notag
\end{equation}
     i.e. unit velocity can be assumed.
\begin{Remark}
Thus we have harmonized the scales along the coordinate axes.
\end{Remark}

     The equation systems of a helix-like geodesic curves $g(x(t),y(t),z(t))$ if $0<|w| <1 $:
    \begin{equation}
    \begin{gathered}
   x(t)=\frac{2c}{w} \sin{\frac{wt}{2}}\cos\Big( \frac{wt}{2}+\alpha \Big),\ \ 
   y(t)=\frac{2c}{w} \sin{\frac{wt}{2}}\sin\Big( \frac{wt}{2}+\alpha \Big), \notag \\
   z(t)=wt\cdot \Big\{1+\frac{c^2}{2w^2} \Big[ \Big(1-\frac{\sin(2wt+2\alpha)-\sin{2\alpha}}{2wt}\Big)+ \\
   +\Big(1-\frac{\sin(2wt)}{wt}\Big)-\Big(1-\frac{\sin(wt+2\alpha)-\sin{2\alpha}}{2wt}\Big)\Big]\Big\} = \\
   =wt\cdot \Big\{1+\frac{c^2}{2w^2} \Big[ \Big(1-\frac{\sin(wt)}{wt}\Big)
   +\Big(\frac{1-\cos(2wt)}{wt}\Big) \sin(wt+2\alpha)\Big]\Big\}. \tag{1.7}
  \end{gathered}
  \end{equation}
In the cases $w=0$ the geodesic curve is the following:
\begin{equation}
x(t)=c\cdot t \cos{\alpha}, \ \ y(t)=c\cdot t \sin{\alpha}, \ \ z(t)=\frac{1}{2} ~ c^2 \cdot t^2 \cos{\alpha} \sin{\alpha}. \tag{1.8}
\end{equation}
The cases $|w|=1$ are trivial: $(x,y)=(0,0), \ z=w \cdot t$.

\begin{Definition}
The distance $d(P_1,P_2)$ between the points $P_1$ and $P_2$ is defined by the arc length of the geodesic curve 
from $P_1$ to $P_2$.
\end{Definition}

 \subsection{On the geodesic ball}
 
 In our work \cite{Sz07-2} we have introduced the followin definitions:
  \begin{Definition}
 The geodesic sphere of radius $R$ with centre at the point $P_1$ is defined as the set of all points 
 $P_2$ in the space with the condition $d(P_1,P_2)=R$. Moreover, we require that the geodesic sphere is a simply connected 
 surface without selfintersection 
 in the $\NIL$ space.
 \begin{Remark}
 We will see that this last condition depends on radius $R$.
 \end{Remark}
 \end{Definition}
 \begin{Definition}
 The body of the geodesic sphere of centre $P_1$ and of radius $R$ in the $\NIL$ space is called geodesic ball, denoted by $B_{P_1}(R)$,
 i.e. $Q \in B_{P_1}(R)$ iff $0 \leq d(P_1,Q) \leq R$.
 \end{Definition}
 \begin{Remark}
 Henceforth, typically we choose the origin as centre of the sphere and its ball, by the homogeneity of 
 $\NIL$. 
\end{Remark}

We apply the quadratic mapping $\mathcal{M}: \NIL \longrightarrow \mathbf{A}^3$ at (1.5) to the geodesic sphere $S$, 
its $\mathcal{M}$-image is denoted by $S'=\mathcal{M}(S)$. 

We choose a point $P(x(R,\theta,\alpha),y(R,\theta,\alpha),z(R,\theta,\alpha))$ lying on a sphere $S$ of radius $R$ with 
centre at the origin. The coordinates of $P$ are given by parameters $(\alpha\in[-\pi, \pi), ~ \theta\in [-\frac{\pi}{2},\frac{\pi}{2}], ~ R>0)$ (see (1.5), (1.10)), 
its $\mathcal{M}$-image is $P'(x'(R,\theta,\alpha),y'(R,\theta,\alpha),z'(R,\theta,\alpha))\in S'$ where
    \begin{equation}
    \begin{gathered}
   x'(R,\theta,\alpha)=\frac{2c}{w} \sin{\frac{wR}{2}}\cos\Big( \frac{wR}{2}+\alpha \Big), \notag \\ 
   y'(R,\theta,\alpha)=\frac{2c}{w} \sin{\frac{wR}{2}}\sin\Big( \frac{wR}{2}+\alpha \Big), \notag \\
   z'(R,\theta,\alpha)=wR+\frac{c^2R}{2w} - \frac{c^2}{2w^2}\sin{wR}, \ \ (\theta\in [-\frac{\pi}{2},\frac{\pi}{2}]\setminus \{0\}), \\
   \text{if} \ \theta=0 \ \text{then} \ x'(R,0,\alpha)=R \cos{\alpha}, \\ y'(R,0,\alpha)=R \sin{\alpha}, \ z'(R,0,\alpha)=0. \tag{1.9}
  \end{gathered}
  \end{equation}
 We can see from the last equations that $(x')^2+(y')^2=\frac{4c^2}{w^2}\sin^2{\frac{wR}{2}}$ and that the $z'$-coordinate does not depend on the parameter 
 $\alpha$, therefore $S'$ can be generated by rotating the following curve about the $z$ axis (lying in the plane $[x,z]$):
     \begin{equation}
     \begin{gathered}
    X(R,\theta)=\frac{2c}{w} \sin{\frac{wR}{2}}=\frac{2\cos{\theta}}{\sin{\theta}} \sin{\frac{R \sin{\theta}}{2}}, \\ 
    Z(R,\theta)=wR+\frac{c^2R}{2w} - \frac{c^2}{2w^2}\sin{wR}= \\ 
    R\sin{\theta}+\frac{R\cos^2{\theta}}{2\sin{\theta}} - \frac{\cos^2{\theta}}{2\sin^2{\theta}}\sin(R\sin{\theta}), \ \ (\theta\in [-\frac{\pi}{2},\frac{\pi}{2}]\setminus\{0\}); \\
    \text{if} \ \theta=0 \ \text{then} \ X(R,0)=R, \ Z(R,0)=0. \tag{1.10}
   \end{gathered}
   \end{equation}
\begin{Remark}
From the definition of the quadratic mapping $\mathcal{M}$ at (1.5) it follows that the cross section 
of the spheres $S$ and $S'$ with the plane $[x,z]$, is the same curve which is 
specified by the parametric equations (1.10).
\end{Remark}
\begin{Remark}
The parametric equations of the geodesic sphere of radius $R$ can be generated from (1.10) by $\NIL$ rotation (see (1.4)).
\end{Remark}

We have denoted by $B(S)$ the body of the $\NIL$ sphere $S$ and by $B(S')$ the body of the sphere $S'$,
furthermore we have denoted their volumes by $Vol(B(S))$ and $Vol(B(S'))$, respectively.

In \cite{Sz07-2} we have proved the the following theorem:
\begin{theorem}
The geodesic sphere and ball of radius $R$ exists in the $\NIL$ space if and only if $R \in [0,2\pi].$
\end{theorem}

We obtain the volume of the geodesic ball of radius $R$ by the following integral (see 1.10):
\begin{equation}
\begin{gathered}
Vol(B(S))=2 \pi \int_0^{\frac{\pi}{2}}X^2 ~ \frac{\mathrm{d}~Z}{\mathrm{d}~\theta} ~ \mathrm{d}~{\theta}= \\ =2 \pi \int_0^{\frac{\pi}{2}}
\Big (\frac{2\cos{\theta}}{\sin{\theta}} \sin{\frac{(R \sin{\theta})}{2}}\Big)^2 \cdot \Big(-\frac{1}{2} \frac{R \cos^3{\theta}}{\sin^2{\theta}}+ 
\frac{\cos{\theta}\sin{(R \sin{\theta})}}{\sin{\theta}}+ \\ +\frac{\cos^3{\theta}\sin{(R \sin{\theta})}}{\sin^3{\theta}}-
\frac{1}{2} \frac{R \cos^3{\theta}\cos{(R \sin{\theta})}}{\sin^2{\theta}}\Big){\mathrm{d}~\theta}. \tag{1.11}
\end{gathered}
\end{equation}

The $\NIL$ sphere of radius $R$ is generated by the $\NIL$ rotation about the axis $z$ 
(see the equation system (1.10) and remarks (1.4), (1.5)). The parametric equation system of the geodesic sphere $S(R)$ in 
our model:
\begin{equation}
\begin{gathered}
   x(R,\theta,\phi)=\frac{2c}{w} \sin{\frac{wR}{2}}\cdot \cos{\phi}= \frac{2\cos{\theta}}{\sin{\theta}} 
   \sin{\frac{R \sin{\theta}}{2}}\cdot \cos{\phi}, \notag \\ 
   y(R,\theta,\phi)=\frac{2c}{w} \sin{\frac{wR}{2}}\cdot \sin{\phi}= \frac{2\cos{\theta}}{\sin{\theta}} 
   \sin{\frac{R \sin{\theta}}{2}}\cdot \sin{\phi}, \notag \\
   z(R,\theta,\phi)=wR+\frac{c^2R}{2w}-\frac{c^2}{2w^2}\sin{wR} +\frac{1}{4}\Big( \frac{2c}{w} 
   \sin{\frac{wR}{2}}\Big)^2 \sin{2\phi}= \\
   =R\sin{\theta}+\frac{R\cos^2{\theta}}{2\sin{\theta}} - \frac{\cos^2{\theta}}{2\sin^2{\theta}}\sin(R\sin{\theta})+
   \frac{1}{4}\Big( \frac{2\cos{\theta}}{\sin{\theta}} \sin{R \frac{\sin{\theta}}{2}}\Big)^2 \sin{2\phi} \notag \\
  \end{gathered}
\end{equation}  
 \begin{equation}
\begin{gathered} 
  -\pi<\phi\leqq \pi, \ \ -\frac{\pi}{2}\leqq \theta \leqq \frac{\pi}{2} \ \text{and} \ \theta \ne 0. \\
   \text{if} \ \theta=0 \ \text{then} \ x(R,0,\phi)=R \cos{\phi}, \ \ y(R,0,\phi)=R \sin{\phi}, \\ 
   \ z(R,0,\phi)=\frac{1}{2} ~ R^2 \cos{\phi} \sin{\phi}. \tag{1.12}
\end{gathered}
\end{equation}
The following theorem was obtained by the derivatives of these parametrically represented functions (by intensive and careful computations with {\it Maple} 
through the second fundamental form) (see \cite{Sz07-2}):
\begin{theorem}
The geodesic $\NIL$ ball $B(S(R))$ is convex in affine-Euclidean sense in our model if and only if $R \in [0,\frac{\pi}{2}]$. 
\end{theorem}
\subsection{The discrete translation group L(Z,~{\it{k}})}
We consider the $\NIL$ translations defined in (1.1) and (1.3) and choose two arbitrary translations 
\begin{equation}
\begin{gathered}
\tau_1=
\begin{pmatrix}
1&t_1^1&t_1^2&t_1^3 \\
0&1&0&0 \\
0&0&1&t_1^1 \\
0&0&0&1 \\
\end{pmatrix} \ \text{and} \ 
\tau_2=
\begin{pmatrix}
1&t_2^1&t_2^2&t_2^3 \\
0&1&0&0 \\
0&0&1&t_2^1 \\
0&0&0&1 \\
\end{pmatrix},
\end{gathered} \tag{1.13}
\end{equation}
now with upper indices for coordinate variables.
We define the translation $(\tau_3)^k,\ \ (k $ $\in $ $\mathbb{N},$ $k $ $\ge 1)$ by the following commutator:
\begin{equation}
(\tau_3)^k=\tau_2^{-1}\tau_1^{-1}\tau_2\tau_1=
\begin{pmatrix}
1 & 0 & 0 & -t_2^1 t_1^2+t_1^1 t_2^2 \\
0&1&0&0\\
0&0&1&0\\
0&0&0&1
\end{pmatrix}. \tag{1.14}
\end{equation}
If we take integers as coefficients, their set is denoted by $\mathbb{Z}$, then we will generate the discrete group $ (\langle \tau_1,\tau_2 \rangle,k) $ 
denoted by 
$\mathbf{L}(\tau_1,\tau_2,k )$ or by $\mathbf{L}(\mathbb{Z},k)$. 

We know that the orbit space $\NIL /\mathbf{L}(\mathbb{Z},k)$ is a compact manifold, i.e. a $\NIL$ space form.
\begin{definition}
The $\NIL$ point lattice $\Gamma_P (\tau_1,\tau_2,k )$ is a discrete orbit of point $P$ in the $\NIL$ space under the group $\mathbf{L}(\tau_1,\tau_2,k )$= 
$\mathbf{L}(\mathbb{Z},k)$ 
with an arbitrary starting point $P$ for all $(k \in \mathbb{N},$ $k \ge 1)$.
\end{definition}
\begin{Remark}
For simplicity, we have chosen the origin as starting point, by the homogeneity of $\NIL$. 
\end{Remark}
\begin{Remark}
We can assume that $t_1^2=0$, i.e. the image of the origin by the translation $\tau_1$ lies on the plane $[x,z]$.
\end{Remark}
In the following we investigate the most important case $k=1$ where $\tau_3$ correspond to the fibre translation, i.e. $\tau_3=\tau_2^{-1}\tau_1^{-1}\tau_2\tau_1$.

We illustrate the action of $\mathbf{L}(\mathbb{Z},1)$ on the $\NIL$ space in Fig.~1. We consider a non-convex polyhedron 
$\mathcal{F}=OT_1T_2T_3T_{12}T_{21}T_{23}T_{213}T_{13}$, in Euclidean sense, which is determined by translations $\tau_1,\tau_2,\tau_3$.
This polyhedron determines a solid $\widetilde{\mathcal{F}}$ in the $\NIL$ space whose images under $\mathbf{L}(\mathbb{Z},1)$ fill the $\NIL$ space just
once, i.e.without gap and overlap. 

{\it Analogously to the Euclidean integer lattice and parallelepiped, the solid $\widetilde{\mathcal{F}}$ can be called 
$\NIL$ parallelepiped.}

$\widetilde{\mathcal{F}}$ is a {\it fundamental domain} of $\mathbf{L}(\mathbb{Z},1)$. 
The homogeneous coordinates of the vertices of $\widetilde{\mathcal{F}}$ can be determined in our affine model by the translations
(1.13) and (1.14) with the parameters $t_i^j,~~k=1 \ (i\in\{1,2\}, \ j \in \{1,2,3\})$ (see Fig.~1 and (1.15)). 
\begin{figure}[ht]
\centering
\includegraphics[width=8cm]{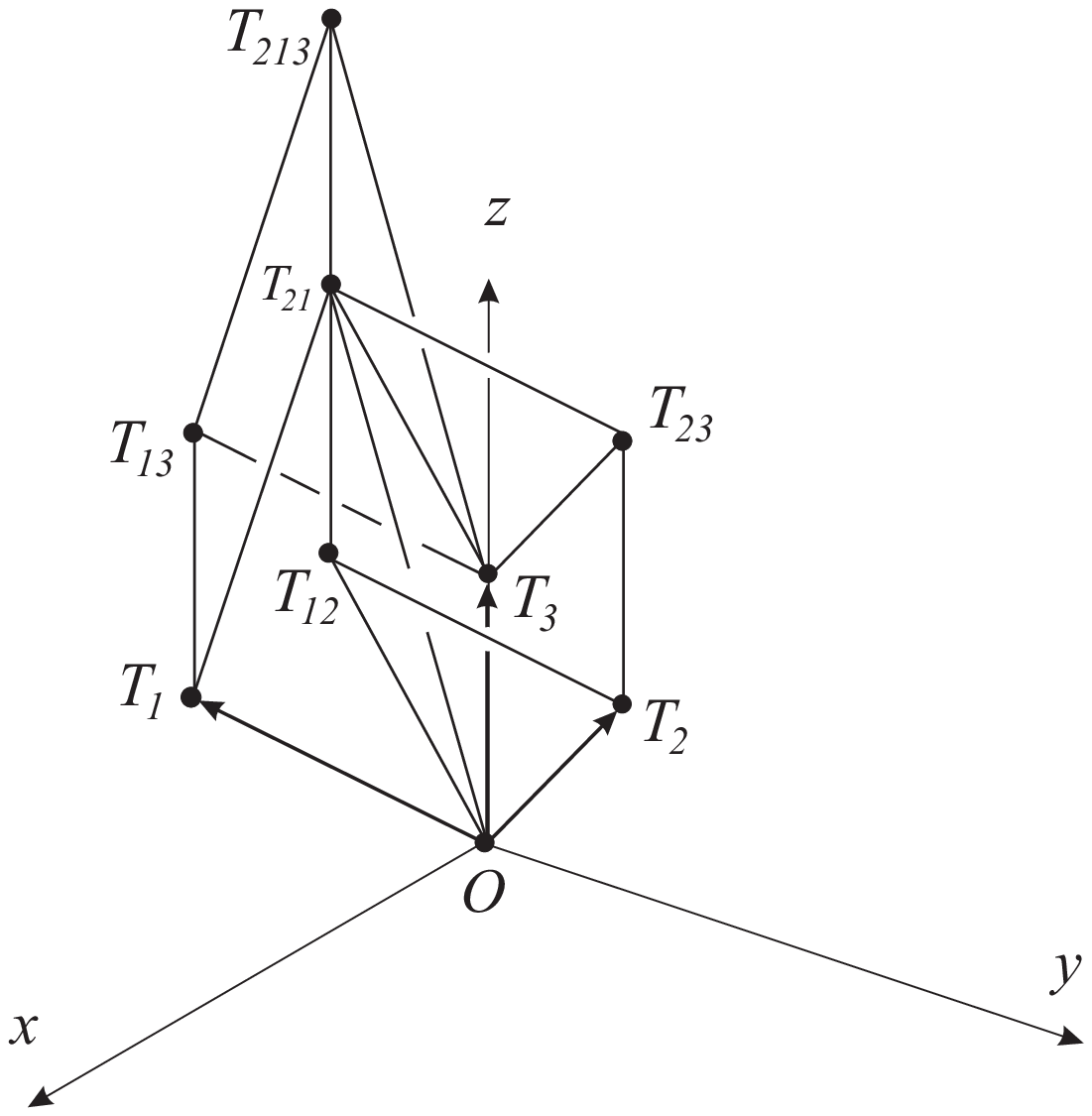}
\caption{}
\label{}
\end{figure}
\begin{equation}
\begin{gathered}
T_1(1,t_1^1,0,t_1^3), \ T_2(1,t_2^1,t_2^2,t_2^3), \ T_3(1,0,0,\frac{t_1^1 t_2^2}{k}), \\ 
T_{13}(1,t_1^1,0,\frac{t_1^1 t_2^2}{k}+t_1^3), \ T_{12}(1,t_1^1+t_2^1,t_2^2,t_2^3+t_1^3), \\ 
T_{21}(1,t_1^1+t_2^1,t_2^2,t_1^1 t_2^2+t_1^3+t_2^3), \ 
T_{23}(1,t_2^1,t_2^2,t_2^3+\frac{t_1^1 t_2^2}{k}), \\ T_{213}=T_{231}(1,t_1^1+t_2^1,t_2^2, (k+1) \frac{t_1^1 t_2^2}{k}+t_1^3+t_2^3). \tag{1.15}
\end{gathered}
\end{equation}
We have determined in \cite{Sz07-2} the volume of $\widetilde{\mathcal{F}}$ by the following integral:
\begin{equation}
Vol(\widetilde{\mathcal{F}})=\int_0^{t_2^2} \int_0^{t_1^1} ~|t_1^1 \cdot t_2^2| ~ \mathrm{d}{x}\mathrm{d}{y}=(t_1^1 \cdot t_2^2)^2. \tag{1.16}
\end{equation}
From this formula it can be seen that the volume of the $\NIL$ parallelepiped depends on two parameters, i.e. on its projection onto the $[x,y]$ plane.
\section{The lattice-like geodesic ball coverings}

A family of subsets $\mathcal{K}=(K_i)_{i\in I}$ of $\NIL$, $I$ is a set of indices is called {\it covering} of $\NIL$ if each point of
$\NIL$ belongs at least one of the set $\{K_i,~ i\in I \}$, i.e. $\NIL = \bigcup_{i \in I} K_i$. A covering of $\NIL$ space is a 
{\it lattice covering} if it is of the form $\mathcal{K}=(K_0+\mathbf{v})_{v \in \mathbf{L}(\mathbb{Z},k)}$ where $K_0$ is a element of 
$\{K_i,~ i\in I \}$  i.e. the lattice coverings those coverings which cover $\NIL$ by translated copies of a single body $K_0$ and in addition 
the translates are vectors of a lattice $ \mathbf{L}(\mathbb{Z})$.

In following, we are only considering lattice coverings consisting of geodesic balls of the $\NIL$.  
Let $\mathcal{B}^c_\Gamma(R)$ denote a geodesic ball covering of $\NIL$ space with balls $B^c(R)$ of radius $R$ where their 
centres give rise to a $\NIL$ point lattice
$\Gamma(\tau_1,\tau_2,1)$. $\widetilde{\mathcal{F}}_0$ is an arbitrary $\NIL$ parallelepiped of this lattice (see (1.13), (1.14)). 
The images of $\widetilde{\mathcal{F}}_0$ by our discrete translation group $\mathbf{L}(\tau_1,\tau_2,1)$
cover the $\NIL$ space without overlap. 
\begin{rmrk}
In the $d$-dimensional Euclidean space $\bE^d, \ (d\ge 1)$ an arbitrary lattice $\Gamma(\tau_1,\tau_2,1)$ under the group $ \mathbf{L}(\mathbb{Z},1)$ 
gives a lattice covering of equal
balls if the radius $R$ of the balls is large enough, but this is not true in the $\NIL$ space, because of a geodesic ball
exists in the $\NIL$ space if and only if $R \in [0,2\pi]$ (see Theorem 1.2).
\end{rmrk}
If we start with a lattice covering $\mathcal{B}^c_\Gamma(R)$ and shrink the balls until they finally do not cover the space any more, then the 
threshold value of the shrinking radius $R$ defines the least dense covering of equal balls to a given lattice $\Gamma(\tau_1,\tau_2,1)$.
The thresfold value $R^c$ is called {\it covering radius} of the point lattice $\Gamma(\tau_1,\tau_2,1)$:
\begin{equation}
R^c:=\min \{R: ~ \text{where} ~  \mathcal{B}^c_\Gamma(R) \ \text{lattice covering by} ~ \Gamma(\tau_1,\tau_2,1) \}. \tag{2.1}
\end{equation}

For the density of the packing it is sufficient to relate the volume of the "minimal covering ball"
to that of the solid $\widetilde{\mathcal{F}}_0$. 

Analogously to the Euclidean case it can be defined the density $\delta(R^c,\tau_1,\tau_2,1)$ of the lattice-like geodesic
ball covering $\mathcal{B}_\Gamma(R^c)$:
\begin{definition}
\begin{equation}
\Delta(R^c,\tau_1,\tau_2,1):=\frac{Vol(\mathcal{B}_\Gamma(R^c))}{Vol(\widetilde{\mathcal{F}})}, \tag{2.2}
\end{equation}
\end{definition}
{\it The main problem is that to which lattice $\mathbf{L}(\tau_1,\tau_2,1)$ belongs the minimal density $\Delta_{opt}$}.
We introduce for the "optimal arrangement" the following denotations:
\begin{equation}
\Delta_{opt}(R^c_{opt},\tau^{c}_1,\tau^{c}_2,1):=
\min \Big\{ \frac{Vol(\mathcal{B}_\Gamma(R^c))}{Vol(\widetilde{\mathcal{F}})} \Big\}. \tag{2.3}
\end{equation}
\begin{rmrk}
The covering radius is the radius of the circumsphere of the lattice's Dirichlet-Voronoi polytope, that is the largest distance
between the midpoint and the vertices of its Dirichlet-Voronoi polytope.
\end{rmrk}
\subsection{The lattice-like ball covering of the lattice $\mathbf{L}_{opt}(\tau_1^{opt},\tau_2^{opt},1)$}
First we consider the ball arrangement $\mathcal{B}_\Gamma^{opt}(R_{opt})$ of the {\it densest lattice-like geodesic ball packing} 
in the $\NIL$ space, given by formulas (2.4), (2.5) (see \cite{Sz07-2}). 
\begin{equation}
\begin{gathered}
(a) \ \ d(O,T_1)=2R=d(T_1,T_3), \\
(b) \ \ d(O,T_2)=2R=d(T_2,T_3), \\
(c) \ \ d(T_1,T_2)=2R, \\
(d) \ \ d(O,T_3)=2R.
\end{gathered} \tag{2.4}
\end{equation}
\begin{equation}
\begin{gathered}
t_1^{1,opt} \approx 1.30633820, \ t_1^{3,opt} =R_{opt}, \ R_{opt} \approx 0,73894461; \\
t_2^{1,opt} \approx 0,65316910, \ t_2^{2,opt} \approx 1,13132206, \ t_2^{3,opt} \approx 1.10841692,  \\
T_1^{opt}=(1,t_1^{1,opt},0,t_1^{3,opt}), \ \ T_2^{opt}=(1,t_2^{1,opt},t_2^{2,opt},t_2^{3,opt}).
\end{gathered} \tag{2.5}
\end{equation}
This packing can be generated by the translations $\mathbf{L}_{opt}(\tau_1^{opt},\tau_2^{opt},1)$
where $\tau_1^{opt}$ and $\tau_2^{opt}$ are given by the coordinates $t_i^{j, opt} \ i=1,2; \ j=1,2,3$ (see (2.5)). Thus we obtain the 
neigbouring balls around an arbitrary ball of the packing $\mathcal{B}_\Gamma^{opt}(R_{opt})$ by the lattice $\Gamma(\tau_1^{opt},\tau_2^{opt},1)$. 
We have ball "columns" in $z$-direction and in regular hexagonal projection onto the $[x,y]$-plane.

The $\NIL$ parallelepiped 
$\widetilde{\mathcal{F}}=OT_1^{opt}T_2^{opt}T_3^{opt}T_{12}^{opt}T_{21}^{opt}T_{23}^{opt}T_{213}^{opt}T_{13}^{opt}$ is a {\it fundamental domain} of $\mathbf{L}(\mathbb{Z},1)$.
The homogeneous coordinates of its vertices are known (see Fig.~1 and (1.15)).

We examine the {\it covering radius} $R^c_p$ to the lattice $\Gamma(\tau_1^{opt},\tau_2^{opt},1)$.
\begin{equation}
R^c_p:=\min \{R: ~ \text{where} ~  \mathcal{B}^c_\Gamma(R) \ \text{lattice covering by} ~ \Gamma(\tau_1^{opt},\tau_2^{opt},1) \}. \notag
\end{equation}

It is sufficient to investigate such ball arrangements $\mathcal{B}^c_{\Gamma^{opt}}(R)$ where the balls cover $\widetilde{\mathcal{F}}$ or
the $\NIL$ solid $\widetilde{\mathcal{P}}=OT_1^{opt}T_{12}^{opt}T_2^{opt}T_3^{opt}T_{13}^{opt}T_{21}^{opt}T_{23}^{opt}$ (see Section 1.2).

From (2.4) and (2.5) follows, that the point sets $\{0, T_1^{opt}, T_2^{opt}, T_3^{opt} \}$,  $\{T_3^{opt}, $ $T_1^{opt}, T_{23}^{opt},$ $T_{13}^{opt} \}$, 
$\{T_3^{opt}, $ $T_1^{opt}, T_{23}^{opt},$ $T_2^{opt} \}$, 
$\{T_{12}^{opt}, T_1^{opt}, T_{23}^{opt}, T_2^{opt} \}$, $\{T_{12}^{opt}, T_1^{opt}, $ $T_{23}^{opt}, $ $T_{13}^{opt} \}$, $\{T_{12}^{opt}, T_{21}^{opt}, 
T_{23}^{opt}, T_{13}^{opt} \}$ are congruent by $\NIL$ isometries.
The radius $R$ of each circumscribed ball to the above point sets can be determined by the following system of equation:
\begin{equation}
\begin{gathered}
d(O,C)=d(C,T_3^{opt})=d(C,T_1^{opt})=d(C,T_2^{opt}),
\end{gathered} \tag{2.6}
\end{equation}
where $C(1,c^1,c^2,c^3)$ is the center of the circumscribed ball of the point set 
$\{0, T_1^{opt}, $ $T_2^{opt}, T_3^{opt} \}$ ($d$ is the $\NIL$ distance, see Definition 1.1):
\begin{equation}
\begin{gathered}
c^1 \approx 0.45981062, \ c^2 \approx  0.26547179, \ c^3 \approx  0.79997799, \\
R \approx 0.90293941.
\end{gathered} \tag{2.7}
\end{equation}
\begin{rmrk}
$C(1,c^1,c^2,c^3)$ is a vertex of the Dirichlet-Voronoi domain of the point O in the $\NIL$ space.
\end{rmrk}
$R\in [0,\frac{\pi}{2}]$ thus by Theorem 1.3 the ball of radius $R$ is convex in affin-Euclidean sense.  We form tetrahedra 
$OT_1^{opt}T_2^{opt}T_3^{opt}$, \ $T_3^{opt}T_1^{opt}T_{23}^{opt}T_{13}^{opt}$, \ $T_3^{opt}T_1^{opt}T_{23}^{opt}T_2^{opt}$, \
\ $T_{12}^{opt}T_1^{opt}T_{23}^{opt}T_2^{opt}$, \ $T_{12}^{opt}T_1^{opt}T_{23}^{opt}T_{13}^{opt}$, \ $T_{12}^{opt}T_{21}^{opt}T_{23}^{opt}T_{13}^{opt}$ in Euclidean sense, which fill the $\NIL$ solid 
$\widetilde{\mathcal{P}}=OT_1^{opt}T_{12}^{opt}T_2^{opt}T_3^{opt}T_{13}^{opt}T_{21}^{opt}T_{23}^{opt}$ just once. Their circumscribed congruent $\NIL$ balls are convex thus they
cover the tetrahedra and so the ball arrangement $\mathcal{B}^c_{\Gamma^{opt}}(R)$ cover the $\NIL$ solid $\widetilde{\mathcal{P}}$. 
Thus the radius $R$ of circumscribed ball 
give us the covering radius $R^c_p$ to the lattice $\mathbf{L}_{opt}(\tau_1^{opt},\tau_2^{opt})$ and we get by (1.11), (1.16) and by the Definition
2.2 the following results:
\begin{equation}
\begin{gathered}
Vol(B(R^c_p)) \approx 3.12538516, \ \ Vol(\widetilde{\mathcal{P}})=Vol(\widetilde{\mathcal{F}}) \approx 2.18415656, \\
\Delta(R^c_p,\tau^{opt}_1,\tau^{opt}_2,1):=\frac{Vol(\mathcal{B}_\Gamma(R^c_p))}{Vol(\widetilde{\mathcal{F}})}
\approx 1.43093459. \tag{2.8}
\end{gathered}
\end{equation}
It follows that 
\begin{equation}
\Delta_{opt}(R^c_{opt},\tau^{c}_1,\tau^{c}_2,1) \le \Delta(R^c_p,\tau^{opt}_1,\tau^{opt}_2,1) \approx 1.43093459. \tag{2.9}
\end{equation}
\begin{rmrk}
The density of the least dense lattice-like ball covering in the the Euclidean space is 
$$
\Delta_{opt}(R^c_{opt},\tau^{c}_1,\tau^{c}_2,1) < \Delta_E=\frac{5 \sqrt{5} \pi}{24} \approx 1.46350307.
$$
\end{rmrk}
\subsection{Upper estimation for the covering radius} 
We consider a arbitrary lattice covering $\mathcal{B}^c_\Gamma(R^c)$ where $\Gamma=\Gamma(\tau_1,\tau_2,1)$ and
\begin{equation}
R^c:=\inf \{R: ~ \text{where} ~  \mathcal{B}^c_\Gamma(R) \ \text{lattice covering by} ~ \Gamma(\tau_1,\tau_2,1) \}. \notag
\end{equation}
The {\it fundamental domain} $\widetilde{\mathcal{F}}$ of the translations group $\mathbf{L}_{opt}(\tau_1,\tau_2,1)$ and its vertices with 
their homogeneous coordinates can be seen in our affine model in (1.15) with parameters $t_i^j, \ i\in\{1,2\}, \ j \in \{1,2,3\}$ (see Fig.~1).
We divide the $\NIL$ solid $\widetilde{\mathcal{P}}=OT_1T_{12}T_2T_3T_{13}T_{21}T_{23}$ (have been derived from the 
fundamental domain $\widetilde{\mathcal{F}}$, see Section 1.2) into $\NIL$ solids  
$OT_1T_2T_3$, \ $T_3T_1T_{23}T_{13}$, \ $T_3T_1T_{23}T_2$, \
\ $T_{12}T_1T_{23}T_2$, \ $T_{12}T_1T_{23}T_{13}$, \ $T_{12}T_{21}T_{23}T_{13}$ which are tetrahedra in terms of Euclidean geometry.
It is clear, that one of them contain the centre of its circumscribed $\NIL$ ball. Suppose now, that this "tetrahedron" is $OT_1T_2T_3$
of which circumscribed ball $B(R)$ centred by $C$ passing through the points $O,T_1,T_2,T_3$. We note here that from conditions of the 
$\NIL$ ball follows that the ball $B(R)$ contain the "Euclidean line segment" $OT_3=\sqrt{Vol(\widetilde{\mathcal{F})}}$ (see Section 1.1-2).
From the Definition 2.2 follows, that
\begin{equation}
\Delta(R^c,\tau_1,\tau_2,1) \ge \Delta(R,\tau_1,\tau_2,1) = 
\frac{Vol(\mathcal{B}_\Gamma(R))}{Vol(\widetilde{\mathcal{F}})}=
\frac{Vol(\mathcal{B}_\Gamma(R))}{OT_3^2}. \tag{2.10}
\end{equation}
\begin{rmrk}
A geodesic curve in $\NIL$ space "parallel to the axis $z$" correspond to an "Euclidean line segment" (see (1.5), (1.7)) with same lenght. 
\end{rmrk}

{\it Forther estimation for the density we need to investigate the upper bound of the length of the line segment $OT_3=\sqrt{Vol(\widetilde{\mathcal{F})}}$ in
$B(R)$.} 

Let $S_*(R)$ be a geodesic sphere of radius $R$ with centre at the origin. 
We apply the quadratic mapping $\mathcal{M}: \NIL \longrightarrow \mathbf{A}^3$ at (1.5) to the geodesic sphere $S_*$, 
its $\mathcal{M}$-image is denoted by $S_*'=\mathcal{M}(S_*)$, moreover 
we have denoted by $B(S_*(R))$ the body of the $\NIL$ sphere $S_*(R)$ and by $B(S_*'(R))$ the body of the sphere $S_*'(R)$, 
\begin{lemma}
The length of the vertical chords of $S_*(R)$ do not change at $\mathcal{M}$.
\end{lemma}
The proof of this lemma follows from the definition of the quadratic mapping $\mathcal{M}$.
\begin{lemma}
$B(S_*'(R))=\mathcal{M}(B(S_*(R)))$ is convex in our model in Euclidean sense if and only if $R \in [0,\pi]$.  
\end{lemma}
{\bf Proof}: \ From the Section 1.1 can be seen that $S_*'(R)$ can be generated by rotating the curve $\mathcal{C}(\theta)=(X(R,\theta),Z(R,\theta))$ 
$(\theta\in [-\frac{\pi}{2},\frac{\pi}{2}]\setminus\{0\})$ (see 1.10) about the $z$ axis (lying in the plane $[x,z]$) thus the convexity of 
the ball $B(S_*'(R))$ follows from the investigation of the derivative 
\begin{equation}
\begin{gathered}
\frac{dZ(R,\theta))}{d\theta}= \\ = \,{\cos \left( {\theta} \right) \Big( \frac { \left( - \left( \cos
 \left( {\theta} \right)  \right) ^{2}R\sin \left( {\theta}
 \right) +2\,\sin \left( \sin \left( {\theta} \right) R \right) 
 \left( \sin \left( {\theta} \right)  \right) ^{2} \right) }{
 2 \left( \sin \left( {\theta} \right)  \right) ^{3}}} - \\
 \,-{\frac {\left(  \left( \cos
 \left( {\theta} \right)  \right) ^{2}\cos \left( \sin \left( {
  \theta} \right) R \right) t\sin \left( {\theta} \right) +2\,
   \left( \cos \left( {\theta} \right)  \right) ^{2}\sin \left( \sin
   \left( {\theta} \right) R \right)  \right) }{2 \left( \sin \left( 
  {\theta} \right)  \right) ^{3}}\Big)}.
 \end{gathered} \tag{2.11}
\end{equation}
From the first (2.11) and second derivatives of $Z(R,\theta))$ follows that if  
$R \in [0, \pi]$ then the ball $B(S_*'(R))$ convex. If $R \in (\pi, 2 \pi]$ then the equation $\frac{dZ(R,\theta))}{d\theta}=0$ 
possesses a solution and the curve $\mathcal{C}$ has an inflection point in the interval 
$\theta \in (0, \frac{\pi}{2})$. In Fig.~2 can be seen the complete curve $\mathcal{C}(\theta)=(X(R,\theta),Z(R,\theta))$ for the parameter
$R=2 \pi$, then the curve possesses at the point $\theta = \frac{\pi}{6}$ the maximum, $Z(2\pi,\frac{\pi}{6})=\frac{5\pi}{2}$.
\begin{figure}[ht]
\centering
\includegraphics[width=7cm]{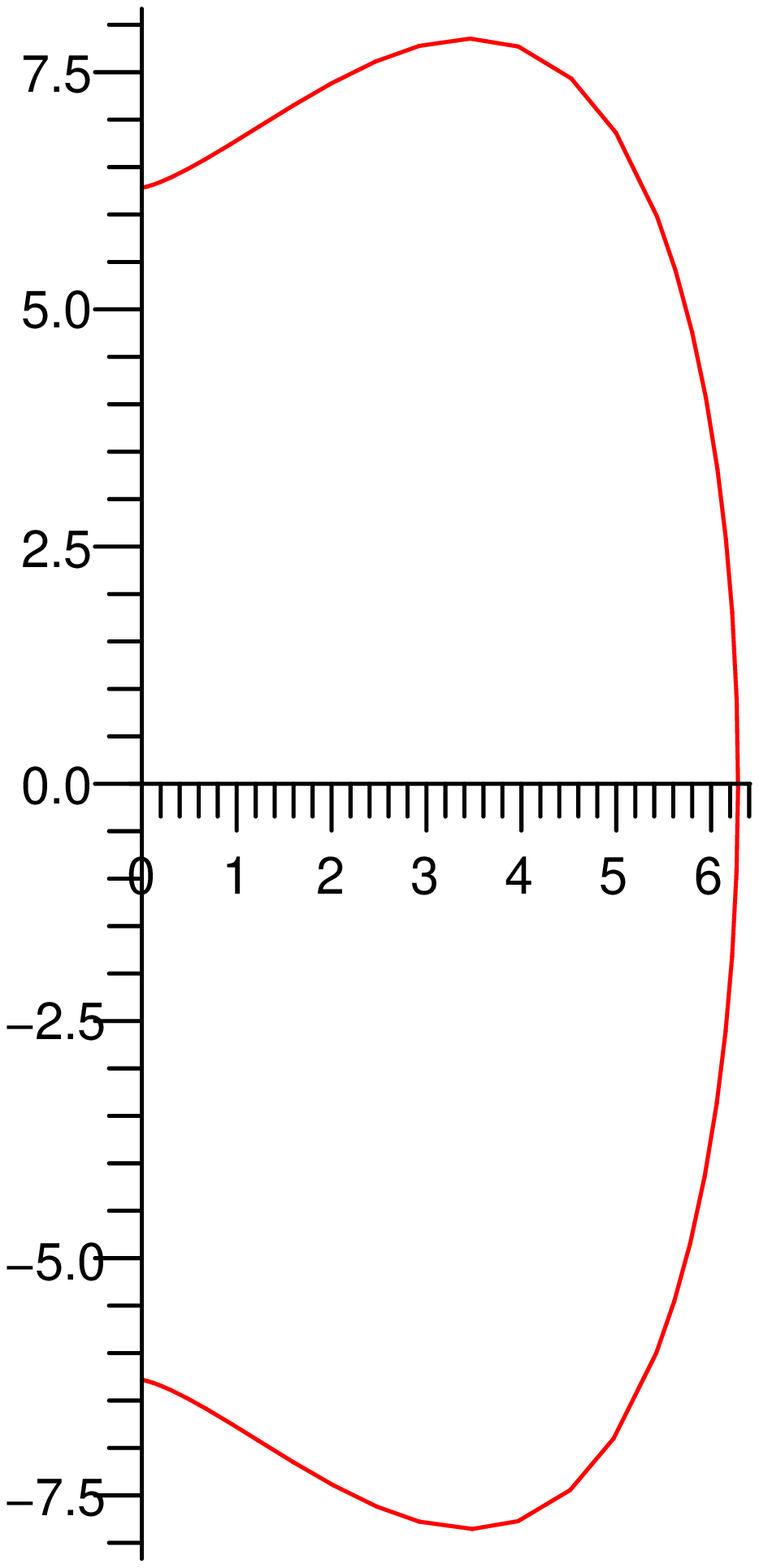}
\caption{}
\label{}
\end{figure}
\begin{lemma}
If $R \in [\frac{\pi}{2},2 \pi]$ then  
\begin{equation}
\begin{gathered}
\Delta(R^c,\tau_1,\tau_2,1) \ge \Delta(R,\tau_1,\tau_2,1) \ge 
\frac{Vol(\mathcal{B}_\Gamma(R))}{(2R)^2} > \\ > \Delta(R^c_p,\tau^{opt}_1,\tau^{opt}_2,1) \approx 1.43093459. \tag{2.12}
\end{gathered}
\end{equation}
\end{lemma}
{\bf Proof}:
\begin{enumerate}
\item $R \in [\frac{\pi}{2}, \pi]$ \newline
We get by (2.10) and by Lemmas 2.6-7 the following inequalities:
\begin{equation}
\Delta(R^c,\tau_1,\tau_2,1) \ge 
\frac{Vol(\mathcal{B}_\Gamma(R))}{OT_3^2} \ge \frac{Vol(\mathcal{B}_\Gamma(R))}{(2R)^2} .
\tag{2.13}
\end{equation}
The function $f(R)=\frac{Vol(\mathcal{B}_\Gamma(R))}{(2R)^2}$ depends only on the parameter $R$ and its graph 
which is increasing on the interval $R \in [\frac{\pi}{2}, \pi]$, can be seen in Fig.~3. Consequently, at the point $R=\frac{\pi}{2}$
the function possesses a minimum, $f(\frac{\pi}{2}) \approx 1.71179510 > \Delta(R^c_p,\tau^{opt}_1,\tau^{opt}_2,1) \approx 1.43093459.$
\begin{figure}[ht]
\centering
\includegraphics[width=7cm]{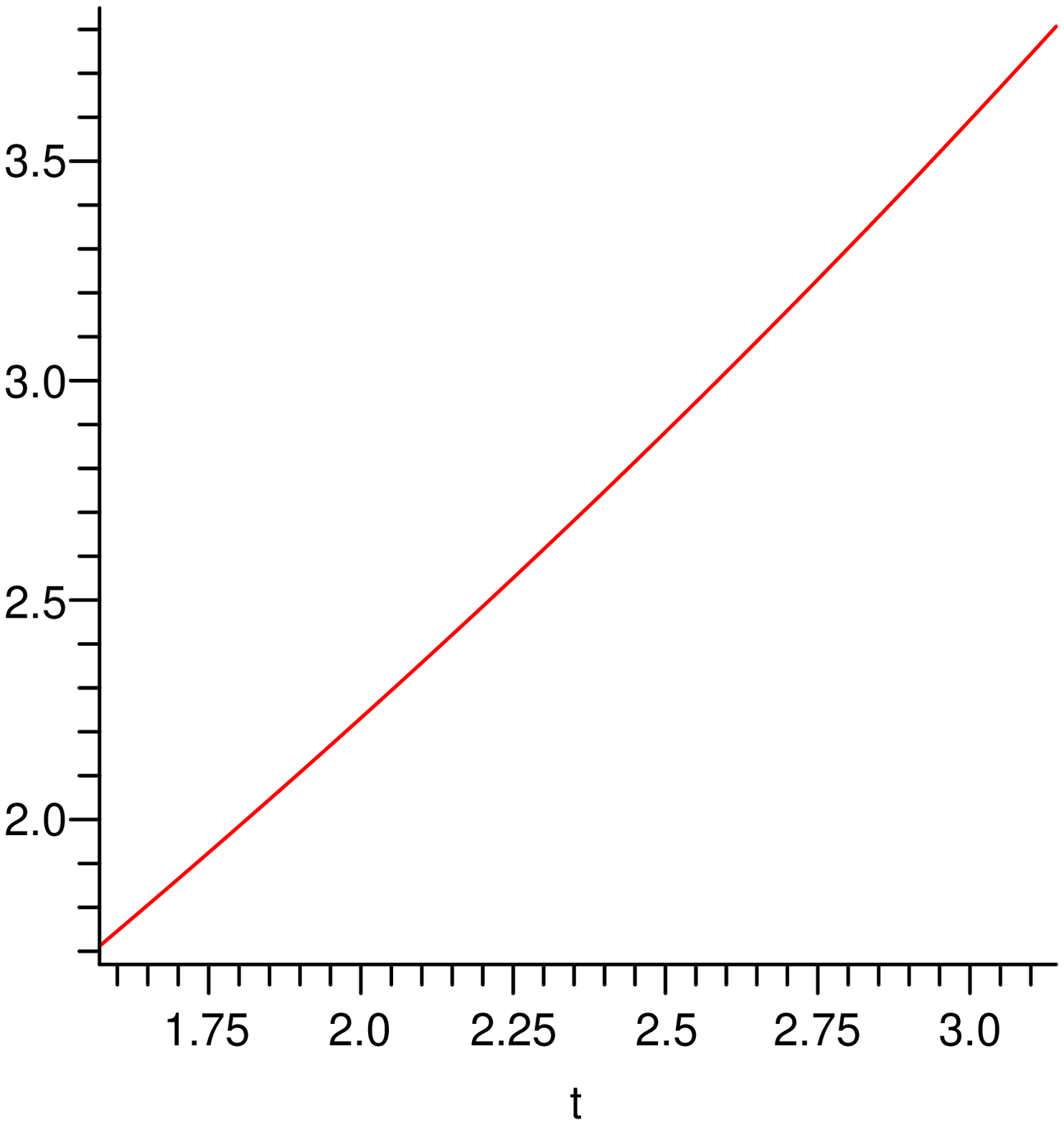}
\caption{}
\label{}
\end{figure}
\item $R \in [\pi, 2 \pi]$ \newline
We divide this interval into two part, $[\pi, 2 \pi]=[\pi, \frac{3 \pi}{2}]\cup (\frac{3 \pi}{2}, 2 \pi]$.
\begin{enumerate}
\item $R \in [\pi, \frac{3 \pi}{2}]$ \newline
In this interval the greatest vertical chord $h_1$ for the balls $B(R)$ is $h_1=\frac{13 \pi}{4} \approx 10.21017613$.
\begin{equation}
\Delta(R^c,\tau_1,\tau_2,1) \ge 
\frac{Vol(\mathcal{B}_\Gamma(R) \cap \widetilde{\mathcal{F}})}{h_1^2}:=f_1(R).
\tag{2.14}
\end{equation}
The Fig.~4.a show its increasing function 
on the interval 
$R \in [\pi, \frac{3 \pi}{2}]$. It is clear that this function  possesses  a  minimum, at the point $R=\pi$. 
$f_1(\pi) \approx 1.441711246$ $> \Delta(R^c_p,\tau^{opt}_1,\tau^{opt}_2,1)$ $ \approx 1.43093459.$ 
\item $R \in (\frac{3 \pi}{2}, 2 \pi]$ \newline
Similarly to 2/a case we examine the possible greatest vertical chord $h_2$ for the balls $B(R)$ on the given interval.
We get that $h_2=5 \pi \approx 15.70796327$ (see Fig.~5).
\begin{equation}
\Delta(R^c,\tau_1,\tau_2,1) \ge 
\frac{Vol(\mathcal{B}_\Gamma(R))}{h_2^2}:=f_2(R).
\tag{2.15}
\end{equation}
The Fig.~4.b shows the graph of the $f_2(R)$ on the interval $[\frac{3 \pi}{2}, 2 \pi]$. It is evident that this function  has a minimum, 
at the point $R=\frac{3 \pi}{2}$. 
$f_2(\frac{3 \pi}{2}) \approx 2.372757787$ $> \Delta(R^c_p,\tau^{opt}_1,\tau^{opt}_2)$ $ \approx 1.43093459.$ 
\end{enumerate}
\end{enumerate}
\begin{figure}[ht]
\centering
\includegraphics[width=5cm]{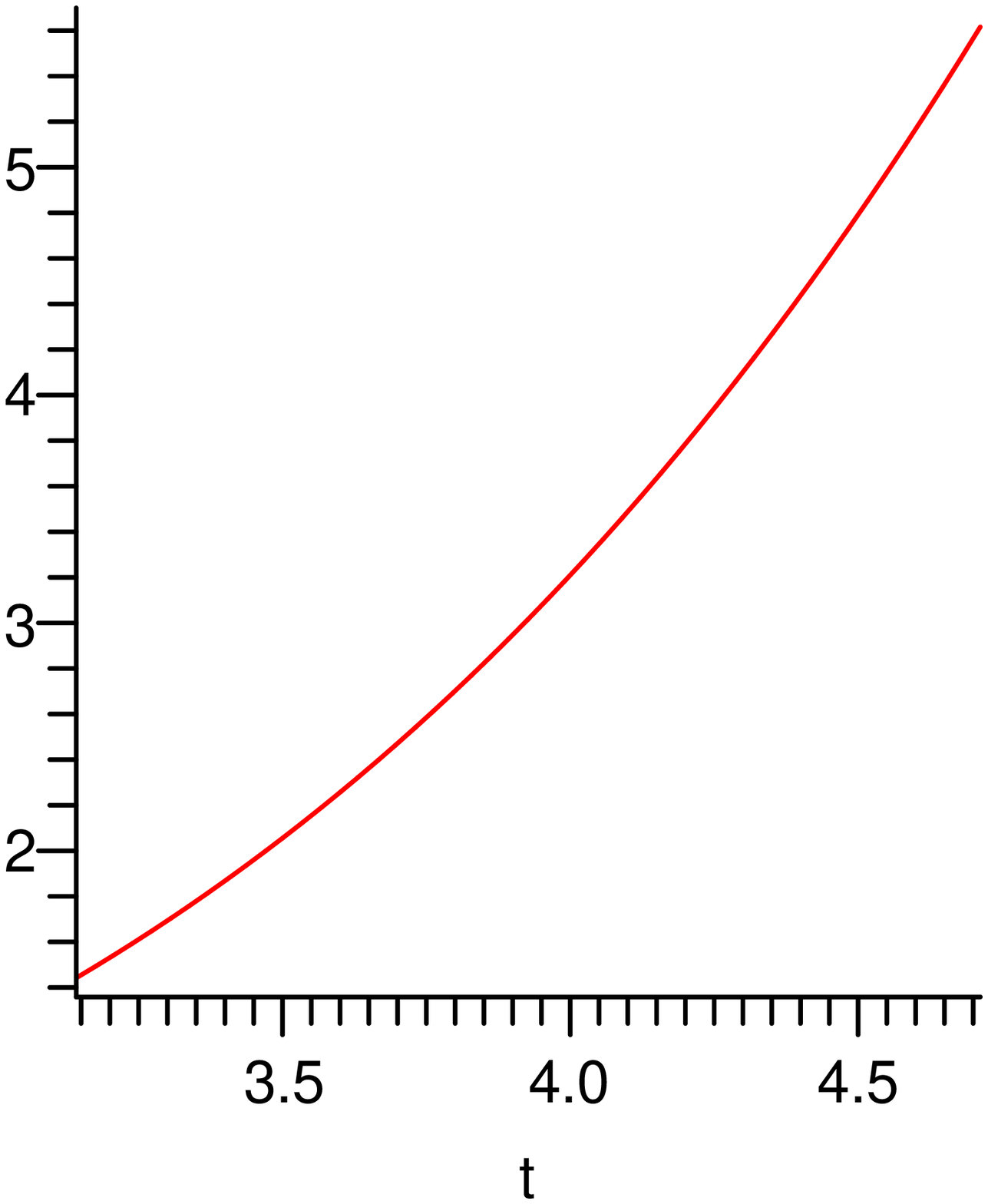} \includegraphics[width=5cm]{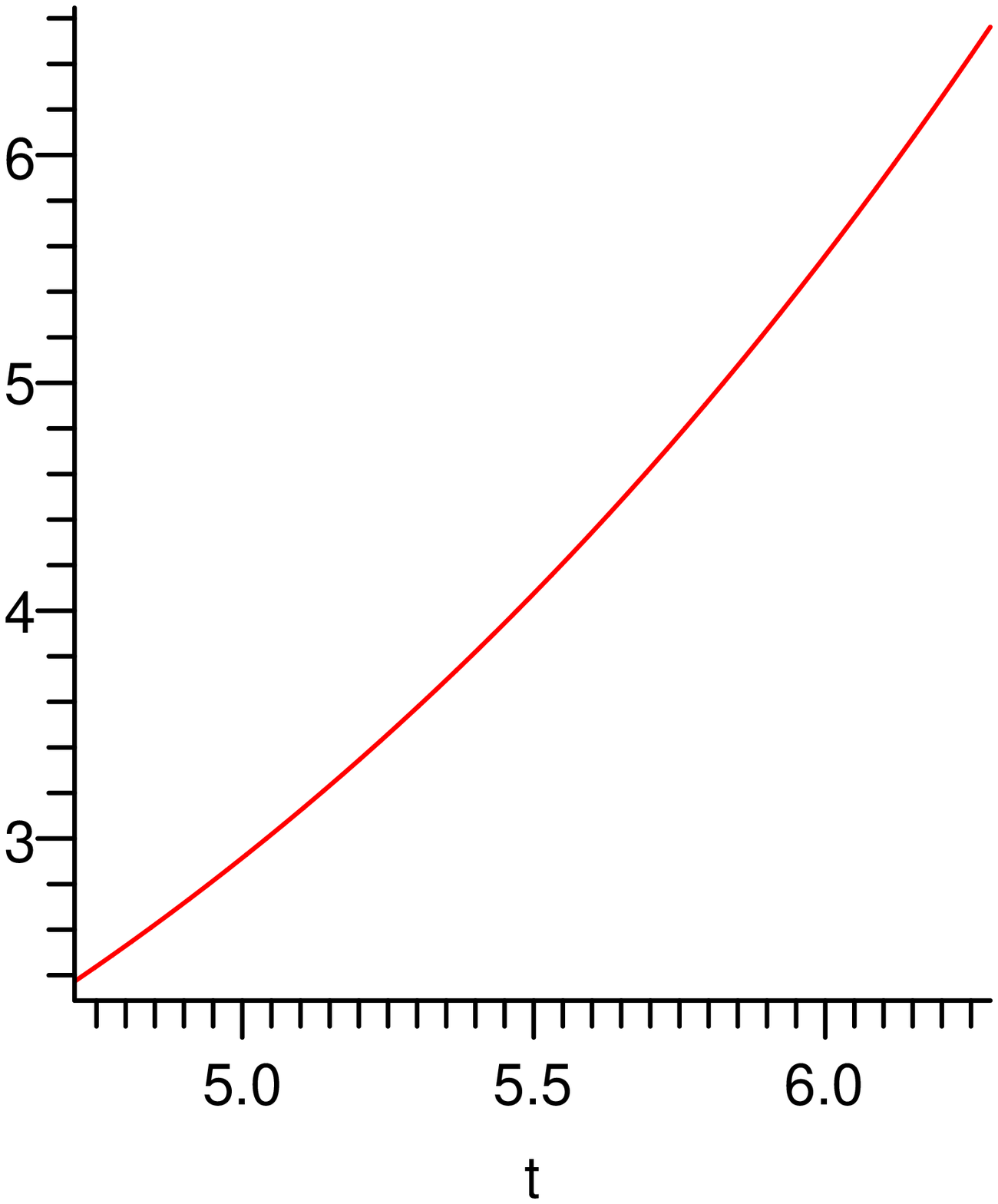}
\caption{a,~b}
\label{}
\end{figure}
Immediate consequence of the Lemma 2.8 is the following Theorem:
\begin{theorem}
The radius $R^c_{opt}$ which belongs to the minimal density $\Delta_{opt}(R^c_{opt},\tau^{c}_1,$ $\tau^{c}_2,1)$
of lattice-like geodesic ball coverings, is found in the interval $[0,\frac{\pi}{2}]$. 
\end{theorem}
\begin{rmrk}
The optimal covering $\NIL$ ball $B(R^c_{opt})$ is convex in Euclidean sense (see Theorem 1.3 and \cite{Sz07-2}).
\end{rmrk}
\subsection{Lower bound to the covering density $\Delta_{opt}(R^c_{opt},\tau^{c}_1,$ $\tau^{c}_2,1)$}
In this section we consider a arbitrary lattice covering $\mathcal{B}^c_\Gamma(R^c)$ where $\Gamma=\Gamma(\tau_1,\tau_2,1)$ (see (2.1-2)).
The {\it fundamental domain} $\widetilde{\mathcal{F}}$ of the translations group $\mathbf{L}(\tau_1,\tau_2,1)$ is given by its vertices 
in our affine model in (1.15) (see Section 1.2 and Fig.~1).
Similarly to Section 2.2 we divide the $\NIL$ solid $\widetilde{\mathcal{P}}=OT_1T_{12}T_2T_3T_{13}T_{21}T_{23}$ (have been derived from the 
fundamental domain $\widetilde{\mathcal{F}}$, see Section 1.2) into $\NIL$ solids  
$OT_1T_2T_3$, \ $T_3T_1T_{23}T_{13}$, \ $T_3T_1T_{23}T_2$, \
\ $T_{12}T_1T_{23}T_2$, \ $T_{12}T_1T_{23}T_{13}$, \ $T_{12}T_{21}T_{23}T_{13}$ which are tetrahedra in Euclidean sense.
It is clear, that one of them contain the centre of its circumscribed $\NIL$ ball and thus can be assumed, that this "tetrahedron" is $OT_1T_2T_3$
of which circumscribed ball $B(R)$ centred by $C(1,c^1,c^2,c^3)$ and passing through the points $O,T_1,T_2,T_3$.  By Theorem 2.10 can be supposed 
that $R \in [0,\frac{\pi}{2}]$ thus $B(R)$ is a convex geodesic $\NIL$ ball in Euclidean sense and contain the Euclidean tetrahedron $OT_1T_2T_3$.

{\it In order to find a lower bound to the covering density we investigate the density function 
\begin{equation}
\Delta(R,\tau_1,\tau_2,1) = 
\frac{Vol(\mathcal{B}_\Gamma(R))}{Vol(\widetilde{\mathcal{F}})}=
\frac{Vol(\mathcal{B}_\Gamma(R))}{OT_3^2} \tag{2.16}
\end{equation}
to a given volume of parallelepiped $Vol(\widetilde{\mathcal{F}})=OT_3^2$. We have to find the minimum radius of the circumscribed ball of the $\NIL$ solid 
$OT_1T_2T_3$ if the Euclidean line segment $OT_3$ is given.} 

We project the points $T_1$ and $T_2$ parallel to the $z$ axis onto the equidistance surface of $O$ and $T_3$ 
which is a hyperbolic paraboloid in our model with equation
$2z-xy=t_1^1 \cdot t_2^2$. Their images are $T_1^p$ and $T_2^p$.   
\begin{lemma}
The radius $R^p$ of the circumscribed ball with centre $C^p(1,c^1_p,c^2_p,c^3_p)$ of the $\NIL$ solid $OT_1^pT_2^pT_3$ which is a tetrahedron in terms of the Euclidean geometry, at most $R$. 
\end{lemma}
{\bf Proof}: From the conditions of the $\NIL$ balls follows (see Section 1.1, \cite{Sz07-2} and Fig.~2) that vertical line segments
from the point upto the equidistance surface are contained by the circumscribed ball $B(R)$, thus $R^p \le R$. Moreover, during this projection
the volume $Vol(\widetilde{\mathcal{F}})=OT_3^2=|t_1^1 \cdot t_2^2|^2$ does not change thus 
\begin{equation}
\Delta(R,\tau_1,\tau_2,1) \ge \frac{Vol(B(R^p))}{Vol(\widetilde{\mathcal{F}})}=
\frac{Vol(B(R^p))}{OT_3^2}. \tag{2.17}
\end{equation}
\begin{rmrk}
Here we do not examine whether $\mathcal{B}_\Gamma(R^p), \ \Gamma=\Gamma(\tau_1^p,\tau_2^p,1)$ is a covering or not.
\end{rmrk}
We consider the Euclidean plane $\alpha$ which passes through the point $T_2^p$ and perpendicular to the axis $y$. If we move $T_2^p$ in this plane then
$t_2^2$ is constant thus $Vol(\widetilde{\mathcal{F}})=OT_3^2=|t_1^1 \cdot t_2^2|^2$ does not change. By Theorem 2.10 and by Remark 2.11 we have obtained the following 
\begin{lemma}
The radius $R^p$ of the circumscribed ball of the $\NIL$ solid $OT_1^pT_2^pT_3$ is minimal at the above moving 
of the point $T_2^p$ if $\alpha$ touch the $\NIL$ ball $B(R^p)$ (see Fig.~5). 
\end{lemma}
\begin{rmrk}
At the "minimal position" of the point $T_2^p$ is $c^1_p=t_2^1$.
\end{rmrk}
\begin{figure}[ht]
\centering
\includegraphics[width=6cm]{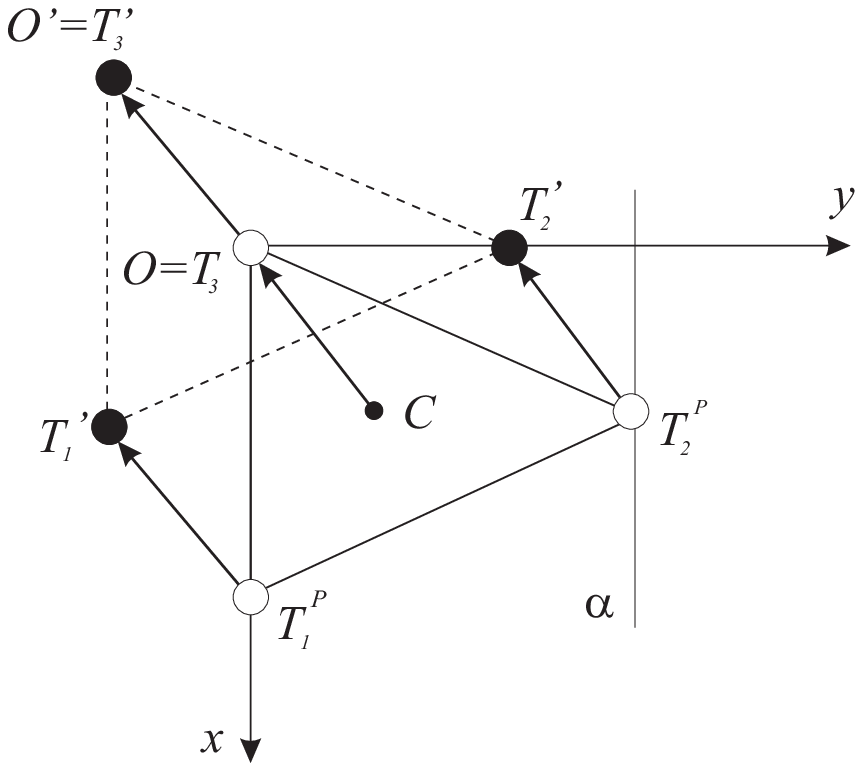}
\caption{}
\label{}
\end{figure}
Forthermore we shall decrease the radius of the circumcribed ball while the volume of the parallelepiped is constant.
\begin{figure}[ht]
\centering
\includegraphics[width=13cm]{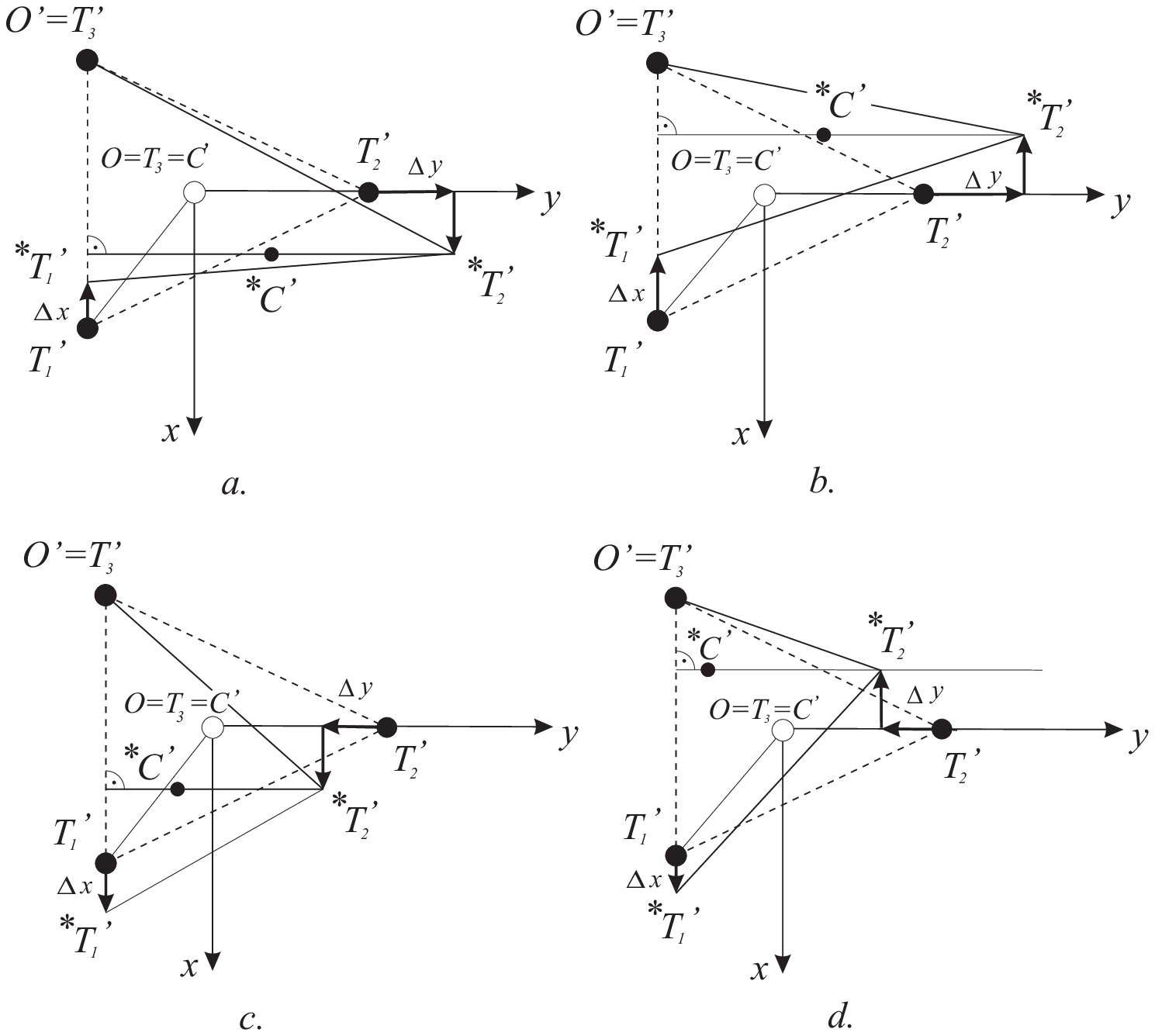}
\caption{}
\label{}
\end{figure}
We translate the $\NIL$ solid $OT_1^pT_2^pT_3$ and its circumscribed ball with $\NIL$ translation $\overrightarrow{C^pO}$, further 
we apply the quadratic mapping $\mathcal{M}: \NIL \longrightarrow \mathbf{A}^3$ at (1.5) to this arrangement (see Fig.~5). 
At these transformations the volume of the parallelepiped $Vol(\widetilde{\mathcal{F}})=OT_3^2$ and the volume of the circumscibed ball 
$Vol(B(R^p))$ do not change.
Thus we get a solid 
$O'T_1'T_2'T_3'$ in affin space $\mathbf{A}^3$ with its circumscribed ball of radius $R^p$ centred by origin. This ball is convex in Euclidean sense
(see Theorem 2.10). Moreover the points $T_1'$, $T_2'$ lie in the plane $[x,y]$ and $OT_3=O'T_3'=\frac{1}{2} \mathcal{A}(H'T_1'T_2')$. Here we have denoted 
by $H' \in [x,y]$ the midpoint of the line segment $O'T_3'$ and by $\mathcal{A}(H'T_1'T_2')$ the area of the Euclidean triangle $H'T_1'T_2'$. 

Working in analogy with what we know from Euclidean geometry, if we fix the volume of the parallelepiped 
$Vol(\widetilde{\mathcal{F}})=OT_3=O'T_3'=\frac{1}{2}$ $\mathcal{A}(H'T_1'T_2')$ then by the Fig.~6. a,b,c,d can be derived the following
\begin{lemma}
The radius $R^p$ of the circumscribed ball of the solid $O'T_1'T_2'T_3'$ to a given volume of parallelepiped is minimal if $H'T_1'=H'T_2'$ where 
$H'T_1'$ and $H'T_2'$ are Euclidean line segments.
\end{lemma}

We have to examine these arrangements. To each $R^p$ can be determined the line segment $O'T_3'$ thus can be investigated the density function 
$\Delta(R^p)$ of these orders and we need to test this function for a possible minimum if $R^p \in [0, \frac{\pi}{2}]$.  
$\Delta(R^p)$ can be examined by careful computation with Maple. The graph of the function $\Delta(R^p)$ can be seen in Fig.~7 and
we get the following resuls:
\begin{equation}
   \begin{gathered}
   \Delta(R^p_{min}) \approx 1.36278112, \ \ R^p_{min} \approx 0.85847445.
   \end{gathered} \tag{2.18}
   \end{equation}
Note that, it is easy to prove, that the ball arrangement belonging to the above given lattice does not yield a geodesic ball covering in the $\NIL$ space. 
\begin{figure}[ht]
\centering
\includegraphics[width=6cm]{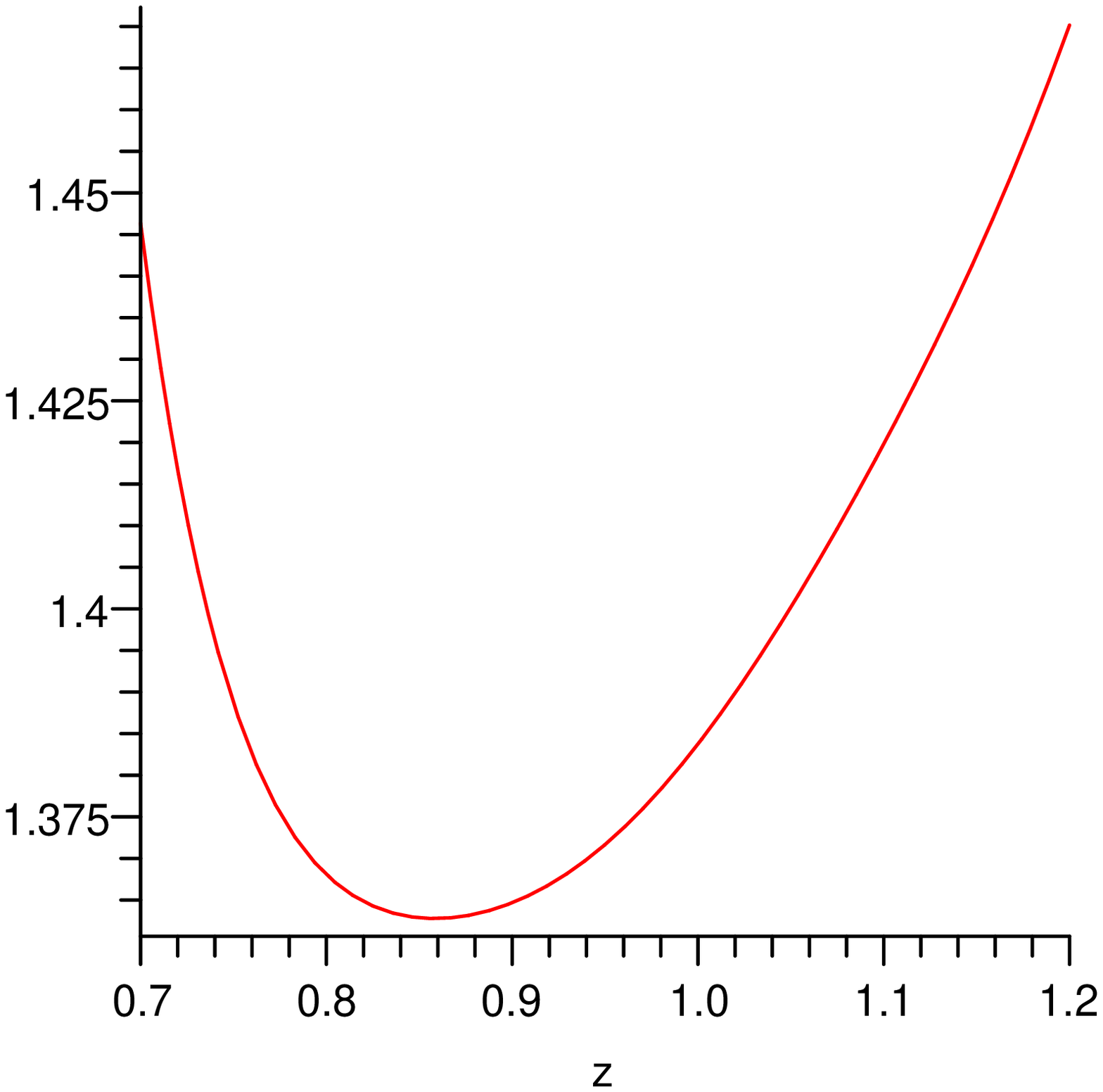}\includegraphics[width=6cm]{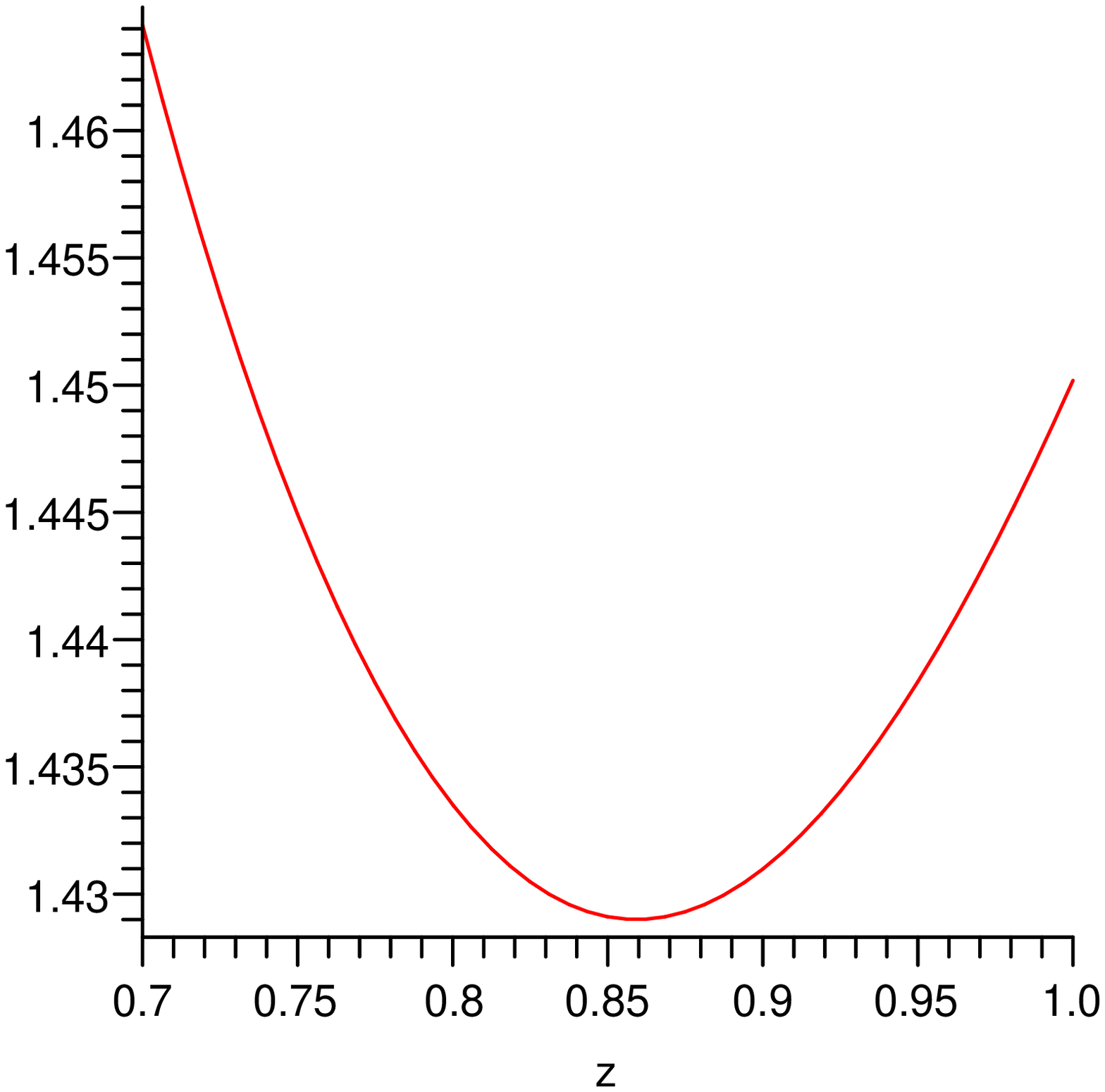}
\caption{a, b}
\label{}
\end{figure}

\begin{corollary}
Consequently, we have obtained a lower bound to the covering density:
\begin{equation}
\Delta_{opt}(R^c_{opt},\tau^{c}_1,\tau^{c}_2,1) > \Delta(R^p_{min}) \approx 1.36278112. \tag{2.19}
\end{equation}
\end{corollary}
\section{Conjecture for the least dense lattice-like \\ ball covering in the $\NIL$ space}
First we consider a $\NIL$ lattice $\Gamma$ which is generated by the translations $\mathbf{L}(\tau_1,$ $\tau_2,1)$
where $\tau_1$ and $\tau_2$ are given by the coordinates $t_i^{j} \ i=1,2; \ j=1,2,3$ (see, Fig.~1). Moreover
the points $T_1$ and $T_2$ lie on the equidistance surface of points $O,~T_3$ and $t_2^1=\frac{t^1_1}{2}, \ \ t_2^2=\frac{\sqrt{3} t^1_1}{2}$,
and one immediate consequence, that $|t_3^3|=|t_1^1 \cdot t_2^1|=\sqrt{Vol(\mathcal{F})}=|(t_1^1)^2 \frac{\sqrt{3}}{2}|$.

We have ball "columns" in $z$-direction and in regular hexagonal projection onto the $[x,y]$ plane.
\begin{rmrk}
The lattice $\Gamma_{opt}$ generated by the translations $\mathbf{L}_{opt}(\tau_1^{opt},\tau_2^{opt},1)$ (see (2.5)) is one of the above
lattices.
\end{rmrk}
The radius $R$ of the circumscribed ball of the $\NIL$ solid $OT_1T_2T_3$ to a given parameter $t_1^1$ can be determined by the following 
system of equation:
\begin{equation}
\begin{gathered}
d(O,C)=d(C,T_3)=d(C,T_1)=d(C,T_2), 
\end{gathered} \tag{3.1}
\end{equation}
where $C(1,c^1,c^2,c^3)$ is the center of the circumscribed ball of the point set $\{0, T_1, $ $T_2, T_3 \}$. 

{\it In order to find the "suspected lower covering density" we investigate the density function} 
\begin{equation}
\Delta(R) = 
\frac{Vol(\mathcal{B}_\Gamma(R))}{Vol(\widetilde{\mathcal{F}})}=
\frac{Vol(\mathcal{B}_\Gamma(R))}{|(t_1^1)^2 \frac{\sqrt{3}}{2}|^2}. \tag{3.2}
\end{equation}
To every $R$ can be determined the parameter $t_1^1$ thus can be examined the density function 
$\Delta(R)$ of these arragement and we need to test this function for a possible minimum if $R \in [0, \frac{\pi}{2}]$.  
$\Delta(R)$ can be investigated by careful computation with Maple. The graph of the function $\Delta(R)$ can be seen in Fig.~7.b and
we get the following resuls:
\begin{equation}
   \begin{gathered}
   \Delta(R^c_{min}) \approx 1.42900615, \ \ R^c_{min} \approx 0.86046718,
   \end{gathered} \tag{3.3}
   \end{equation}
\begin{equation}
\begin{gathered}
t_1^{1,min} \approx 1.26001585, \ t_1^{3,min} \approx 0.68746826; \\
t_2^{1,min} \approx 0,63000792, \ t_2^{2,min} \approx 1,09120574, \ t_2^{3,min} \approx 1.03120239,  \\
T_1^{min}=(1,t_1^{1,min},0,t_1^{3,min}), \ \ T_2=(1,t_2^{1,min},t_2^{2,min},t_2^{3,min}).
\end{gathered} \tag{3.4}
\end{equation}
Similarly to the Section 2.2 it is easy to see, that the ball arrangement belonging to the above given lattice 
is a geodesic ball covering in the $\NIL$ space, thus we get the following 
\begin{theorem}
\begin{equation}
1.36278112 \approx \Delta(R^p_{min})< \Delta_{opt}(R^c_{opt},\tau^{c}_1,\tau^{c}_2,1) \le \Delta(R^c_{min}) \approx 1.42900615. \notag
\end{equation}
\end{theorem}

\begin{conjecture}
The least dense lattice-like ball covering in the $\NIL$ space is derived by the
$\NIL$ lattice $\Gamma_{min}$ which is generated by the translations $\mathbf{L}_{min}(\tau_1^{min},$ $\tau_2^{min},1)$
where $\tau_1^{min}$ and $\tau_2^{min}$ are given by the coordinates $t_i^{j,min} \ i=1,2; \ j=1,2,3$.
The minimal covering radius $R^c_{opt}=R^c_{min} \approx 0.86046718$ and 
$$\Delta_{opt}(R^c_{opt},\tau^{c}_1,\tau^{c}_2,1)=\Delta(R^c_{min}) \approx 1.42900615.$$
\end{conjecture}

Our projective method gives us a way of investigation the $\NIL$ space, which suits to study and solve similar problems (see \cite{Sz07-2}).
In this paper we have examined only some problems, but analogous questions in $\NIL$ geometry 
or, in general, in other homogeneous Thurston geometries are timely (see \cite{Sz10-1}, \cite{Sz10-2}, \cite{Sz10-3}, \cite{Sz11}). 

{\bf{Acknowledgement:}}
I thank Prof. Emil Moln\'ar for helpful comments to this paper.


\begin{thebibliography}{MPSz98}
%
 \bibitem[1]{M97}
  {Moln{\'a}r,~E.}
  The projective interpretation of the eight 3-di\-men\-sional homogeneous geometries. 
  \emph{Beitr{\"a}ge zur Algebra und Geometrie (Contributions to Algebra and Geometry),}
 {\bf38} (1997) No.~2, 261--288.
 %
 \bibitem[2]{M06}
 {Moln{\'a}r,~E.}
  On projective models of Thurston geometries, some relevant notes on $\NIL$ orbifolds and manifolds.  
  \emph{Siberian Electronic Mathematical Reports, http://~semr.math.nsc.ru}  
 {\bf7} (2010), 491--498.
 %
 \bibitem[3]{MPSz}
 {Moln{\'a}r,~E.~---~Prok,~I.~---~Szirmai,~J.}
 Classification of tile-transitive 3-simplex tilings and their realizations in homogeneous spaces.
 \emph{Non-Euclidean Geometries, J\'anos Bolyai Memorial Volume}~ Ed. {\sc{Prekopa,~A.}} and {\sc{Moln\'ar, E.}} \
 Mathematics and Its Applications {\bf 581}, Springer (2006), 321--363.
 %
 \bibitem[4]{S}
 {Scott,~P.}
  The geometries of 3-manifolds.  
  \emph{Bull. London Math. Soc.}, {\bf15} (1983) 401--487. (Russian translation: Moscow "Mir" 1986.)
 %
 \bibitem[5]{Sz07-1}
 {Szirmai,~J.}
  The optimal ball and horoball packings to the Coxeter honeycombs in the hyperbolic d-space. 
  \emph{Beitr{\"a}ge zur Algebra und Geometrie (Contributions to Algebra and Geometry),}
  {\bf48} No.~1, 35--47,(2007).
 %
 \bibitem[6]{Sz07-2}
 {Szirmai,~J.}
 The densest geodesic ball packing by a type of $\NIL$ lattices. 
 \emph{Beitr{\"a}ge zur Algebra und Geometrie (Contributions to Algebra and Geometry),}
 {\bf48} No.~2, 383--398, (2007).
 %
 \bibitem[7]{Sz10-1}
 {Szirmai,~J.}
 The densest translation ball packing by fundamental lattices in $\SOL$ space. 
 \emph{Beitr{\"a}ge zur Algebra und Geometrie (Contributions to Algebra and Geometry),}
 {\bf 51(2)} (2010), 353--373. 
 %
 \bibitem[8]{Sz10-2}
 {Szirmai,~J.}
 Geodesic ball packing in $\SXR$ space for generalized Coxeter space groups.
 \emph{Beitr{\"a}ge zur Algebra und Geometrie (Contributions to Algebra and Geometry), to appear 2011}.
 %
 \bibitem[9]{Sz10-3}
 {Szirmai,~J.}
 Geodesic ball packing in $\HXR$ space for generalized Coxeter space groups.
 \emph{Mathematical Communications, to appear 2011}.
 %
 \bibitem[10]{Sz11}
 {Szirmai,~J.}
 Lattice-like translation ball packings in $\NIL$ space, 
 \emph{Manuscript to Publicationes Mathematicae, Debrecen 2011.}
  %
 \bibitem[11]{T}
 {Thurston,~W.~P.} (and {\sc Levy,~S.} editor) 
 \emph{Three-Dimensional Geometry and Topology.}  
 {Princeton University Press,} Princeton, New Jersey, Vol {\bf 1} (1997).
 %
 \end{thebibliography}
\end{document}